02/25/2008

# EVALUATION OF SOME SIMPLE EULER-TYPE SERIES

Khristo N. Boyadzhiev

Department of Mathematics, Ohio Northern University
Ada, Ohio 45810

**ABSTRACT**

We evaluate five series, three of which involve harmonic numbers and one involves Stirling numbers of the first kind. The evaluation of these series is reduced to the evaluation of certain integrals, including the moments of the polylogarithm.

## 1. INTRODUCTION

In this note we evaluate in closed form the Euler-type sums,

$$\sigma_n(p) \equiv n! \sum_{k=1}^{\infty} \frac{H_k^{(p)}}{k\,(k+1)\ldots(k+n)}\,, \tag{1.1}$$

$$\delta_n(p) \equiv n! \sum_{k=1}^{\infty} \frac{1}{k^p\,(k+1)\ldots(k+n)}\,, \tag{1.2}$$

$$\chi\,(p;n,m) \equiv \sum_{k=1}^{\infty} \frac{H_k^{(p)}}{(k+n)\,(k+m)}\,, \tag{1.3}$$

where

$$H_k^{(p)} = 1 + \frac{1}{2^p} + \ldots + \frac{1}{k^p} \tag{1.4}$$

are the generalized harmonic numbers with $H_k^{(1)} = H_k$ , and in (1.3) $m > n \geq 0$ are integers.

We also evaluate the series

$$\tau_n(p) \equiv n! \sum_{k=p}^{\infty} \frac{W_p(k)}{k(k+1)...(k+n)}, \tag{1.5}$$

where $W_p(k)$ is a p – fold convolution of harmonic numbers. The fifth series is

$$\rho_n(m) \equiv n! \sum_{k=m}^{\infty} \frac{1}{k(k+1)...(k+n)} \begin{bmatrix} k \\ m \end{bmatrix} \frac{1}{k!}, \tag{1.6}$$

where $\begin{bmatrix} k \\ m \end{bmatrix}$ are the Stirling (cycle) numbers of the first kind as defined in [4]. This extends the interesting representation of the Riemann zeta function

$$\zeta(m+1) = \sum_{k=m}^{\infty} \begin{bmatrix} k \\ m \end{bmatrix} \frac{1}{k!k} . \tag{1.7}$$

The evaluation of the first three series is based on the evaluation of the integrals

$$\mu(p;r) \equiv \int_0^1 x^{r-1} \mathrm{Li_p}(x) dx, \tag{1.8}$$

where $\mathrm{Li_p}(x)$ is the polylogarithm

$$\mathrm{Li_p}(x) = \sum_{k=1}^{\infty} \frac{x^n}{n^p}, \tag{1.9}$$

with $\mathrm{Li_p}(1) = \zeta(p)$ when $p > 1$ (see [7, p.106]).

Using the identity

$$\frac{n!}{k(k+1)...(k+n)} = \frac{1}{k}\binom{k+n}{n}^{-1} = \sum_{m=0}^{n} \binom{n}{m} \frac{(-1)^m}{k+m} \tag{1.10}$$

(see [4, 5.41], [6, p. 611] ), the series (1.1), (1.5) and (1.6) can be written in an obvious alternative form.

The motivation for this work came from two problems by Ovidiu Furdiu - problem 854 in the College J. Math., 38 (3), 2007, and problem H-653 in the Fibonacci Quarterly, 45 (1), 2007.

## 2. EVALUATION OF THE MOMENTS OF THE POLYLOGARITHM

The following is true.

**Lemma 2.1.** *For every positive integer* $p$ *and every* $r > 0$

$$\mu(p;r) = \int_0^1 x^{r-1} \mathrm{Li}_p(x)dx = \sum_{n=1}^{\infty} \frac{1}{n^p(n+r)} \tag{2.1}$$

$$= \sum_{k=1}^{p-1} \frac{(-1)^{k-1}\zeta(p-k+1)}{r^k} + \frac{(-1)^{p-1}}{r^p}(\psi(r+1)+\gamma) .$$

Here $\zeta(s)$ is the Riemann zeta function, $\psi(s)$ is the digamma function, and $\gamma = -\psi(1)$ is Euler's constant.

*Proof.* The first equality comes from term-wise integration using the series in (1.9). Then we write

$$\mu(p;r) = \sum_{n=1}^{\infty} \frac{r+n-n}{n^p r(n+r)} = \sum_{n=1}^{\infty} \frac{1}{n^p r} - \sum_{n=1}^{\infty} \frac{1}{n^{p-1} r(n+r)} = \frac{1}{r}[\zeta(p) - \mu(p-1;r)] ,$$

and thus $\mu(p;r)$ is expressed in terms of $\mu(p-1;r)$. Repeating the procedure $p-2$ times brings to (2.1) in view of the fact that

$$\mu(1;r) = \sum_{n=1}^{\infty} \frac{1}{n(n+r)} = \frac{1}{r}(\psi(r+1)+\gamma), \tag{2.2}$$

(see [7, p. 14]). When $r = m$ is a positive integer,

$$\psi(m+1) + \gamma = H_m, \tag{2.3}$$

and therefore,

$$\int_0^1 x^{m-1} \mathrm{Li}_p(x)dx = \sum_{k=1}^{p-1} \frac{(-1)^{k-1}\zeta(p-k+1)}{m^k} + \frac{(-1)^{p-1}}{m^p} H_m . \tag{2.4}$$

In particular, for $m = 1$,

$$\int_0^1 \mathrm{Li}_p(x)dx = \sum_{k=1}^{p-1} (-1)^{k-1} \zeta(p-k+1) + (-1)^{p-1}. \tag{2.5}$$

(This computation is from [2]. The result also appeared in the recent online publication [3, 4.4.100p] together with some other interesting integrals involving polylogarithms.)

The integrals $\mu(p;r)$ can be used to evaluate other integrals of similar type. For example, we have the immediate corollary.

**Corollary 2.2.** *For any positive integer* $m$,

$$\int_0^1 (1-x)^m \mathrm{Li}_p(x)dx = \sum_{j=0}^{m} \binom{m}{j} (-1)^j \mu(p;j+1). \tag{2.6}$$

## 3. EVALUATION OF THE FIRST THREE SERIES

The main tool in the evaluation of the series is Euler's integral

$$\frac{n!}{k(k+1)...(k+n)} = \int_0^1 x^{k-1}(1-x)^n dx \quad , \tag{3.1}$$

(see Euler's Beta function in [6], [7]). We also need the generating function of the generalized harmonic numbers

$$\frac{\mathrm{Li}_p(x)}{1-x} = \sum_{k=1}^{\infty} H_n^{(p)} x^n, \tag{3.2}$$

( $|x|<1$). When $p = 1$ we have $\mathrm{Li}_1(x) = -\ln(1-x)$ and therefore,

$$\frac{-\ln(1-x)}{1-x} = \sum_{n=1}^{\infty} H_n x^n . \tag{3.3}$$

**Proposition 3.1.** *For every two integers* $n \geq 2$ *and* $p > 0$,

$$n! \sum_{k=1}^{\infty} \frac{H_k^{(p)}}{k(k+1)...(k+n)} = (n-1) \sum_{j=0}^{n-2} \binom{n-2}{j} (-1)^j \mu(p+1;j+1) . \tag{3.4}$$

*For instance, when* $n = 2$,

$$\sum_{k=1}^{\infty} \frac{H_k^{(p)}}{k(k+1)(k+2)} = \frac{1}{2}\mu(p+1;1) = \frac{1}{2}\left(\sum_{k=1}^{p} (-1)^{k-1}\zeta(p-k+2) + (-1)^p\right)$$

(For $n = 1$ see (3.11) below).

*Proof.* Using (3.1), (3.2), and the fact that

$$\frac{\mathrm{Li}_p(x)}{x}dx = d\,\mathrm{Li}_{p+1}(x), \tag{3.5}$$

we compute

$$\sigma_n(p) = \sum_{k=1}^{\infty} H_n^{(p)} \int_0^1 x^{k-1}(1-x)^n dx = \int_0^1 \frac{\mathrm{Li}_p(x)}{x(1-x)}(1-x)^n dx$$

$$= \int_0^1 (1-x)^{n-1} d\,\mathrm{Li}_{p+1}(x) = (1-x)^{n-1}\,\mathrm{Li}_{p+1}(x)\Big|_0^1 + (n-1)\int_0^1 (1-x)^{n-2}\,\mathrm{Li}_{p+1}(x)\,dx$$

$$= (n-1)\int_0^1 (1-x)^{n-2}\,\mathrm{Li}_{p+1}(x)\,dx, \tag{3.6}$$

and the proposition now follows from (2.6).

Next we evaluate the series $\delta_n(p)$ in (1.2).

**Proposition 3.2.** *For every* $n \geq 1$, $p \geq 2$

$$n!\sum_{k=1}^{\infty} \frac{1}{k^p(k+1)\ldots(k+n)} = n\sum_{j=0}^{n-1} \binom{n-1}{j} (-1)^j \mu(p;j+1) \tag{3.7}$$

*Proof.* We have according to (3.1),

$$\frac{(n-1)!}{(k+1)\ldots(k+n)} = \int_0^1 x^k(1-x)^{n-1}dx ,$$

hence

$$\frac{1}{n}\delta_n(p) = \sum_{k=1}^{\infty} \frac{1}{k^p}\int_0^1 x^k (1-x)^{n-1}\,dx = \int_0^1 \mathrm{Li_p}(x)\,(1-x)^{n-1}\,dx \tag{3.8}$$

and again we refer to (2.6) in order to complete the proof.

Observation. Comparing the last integrals in (3.6) and (3.8) we notice that

$$\sigma_n(p) = \delta_{n-1}(p+1). \tag{3.9}$$

In the next section we describe a procedure for recursive evaluation of the sums $\delta_n(p)$ for all $n \geq 1, p \geq 1$. According to (3.9), this procedure applies also to the sums $\sigma_n$ .

**Proposition 3.3** *For any three positive integers* $p$ , $m > n$ *one has* $\chi(p; n, m) =$

$$\sum_{k=1}^{\infty} \frac{H_k^{(p)}}{(k+n)(k+m)} = \zeta(p+1) - \frac{1}{m-n}\sum_{j=n}^{m-1} j\,\mu(p+1; j). \tag{3.10}$$

When $n = 0$ , $m = 1$ the sum on the right hand side is missing. Thus

$$\sum_{k=1}^{\infty} \frac{H_k^{(p)}}{k(k+1)} = \zeta(p+1)\,, \tag{3.11}$$

$$\sum_{k=1}^{\infty} \frac{H_k^{(p)}}{k(k+2)} = \zeta(p+1) - \frac{1}{2}\mu(p+1; 1) \tag{3.12}$$

$$= \frac{1}{2}\left[\,\zeta(p+1) + \zeta(p) - \zeta(p-1) \ldots + (-1)^p \zeta(2) + (-1)^{p+1}\,\right].$$

*Proof.* Write first

$$\frac{1}{(k+n)(k+m)} = \frac{1}{m-n}\left(\frac{1}{k+n} - \frac{1}{k+m}\right)$$

$$= \frac{1}{m-n}\int_0^1 x^{k+n-1} - x^{k+m-1}\,dx = \frac{1}{m-n}\int_0^1 x^{k+n-1}(1 - x^{m-n})\,dx\,.$$

Then

$$\chi(p;n,m) = \frac{1}{m-n}\int_0^1 \frac{\mathrm{Li_p}(x)}{x(1-x)}(1-x^{m-n})x^n\,dx$$

$$= \frac{1}{m-n}\int_0^1 (1+x+x^2+\ldots+x^{m-n-1})x^n\,d\,\mathrm{Li_{p+1}}(x)$$

$$= \frac{x^n}{m-n}(1+x+\ldots+x^{m-n-1})\,\mathrm{Li_{p+1}}(x)\Big|_0^1 - \frac{1}{m-n}\int_0^1 \mathrm{Li_{p+1}}(x)\,d(x^n+x^{n+1}+\ldots+x^{m-1})$$

$$= \zeta(p+1) - \frac{1}{m-n}\int_0^1 \mathrm{Li_{p+1}}(x)\,(nx^{n-1}+(n+1)x^n+\ldots+(m-1)x^{m-2})dx\,.$$

and the result follows from (2.1).

**Remark 3.4.** The Euler sum

$$\sum_{k=1}^{\infty}\frac{H_k^{(p)}}{k^2} = S(p;2) \tag{3.13}$$

has been studied by many authors (starting with Euler himself) and its vales are known only for $p = 2, 4$ or $p$ odd (see [1]) . It is interesting that the similar sum (3.11), can be evaluated easily for every $p$ .

## 4. A RECURRENCE FOR THE SECOND SERIES (1.2)

**Lemma 4.1** *For all* $n \geq 1$ , $p \geq 1$ *we have the recursive relation*

$$\delta_n(p) = \delta_{n-1}(p) - \frac{1}{n}\delta_n(p-1), \tag{4.1}$$

*whenever the three sums are convergent.*

To see this, write $n! = (n-1)!\,n$ and $n = n+k-k$ to obtain

$$\frac{n!}{k^p(k+1)...(k+n)} = \frac{(n-1)!}{k^p(k+1)...(k+n-1)} - \frac{(n-1)!}{k^{p-1}(k+1)...(k+n)} ,$$

and the lemma follows immediately.

We shall evaluate now $\delta_n(p)$ for $p = 1, 2, 3$ .

Let first $p = 1$ . From (3.1) we find

$$\delta_n(1) = n!\sum_{k=1}^{\infty} \frac{1}{k(k+1)...(k+n)} = \int_0^1 (1-x)^{n-1}\, dx = \frac{1}{n} . \tag{4.2}$$

Next, let $p = 2$ . Then from (4.1)

$$\delta_n(2) = \delta_{n-1}(2) - \frac{1}{n^2} = \ldots = \delta_0(2) - (1 + \frac{1}{2^2} + \ldots + \frac{1}{n^2}) ,$$

and since

$$\delta_0(2) = \sum_{k=1}^{\infty} \frac{1}{k^2} = \zeta(2), \tag{4.3}$$

we have

$$\delta_n(2) = \zeta(2) - H_n^{(2)} . \tag{4.4}$$

Now let $p = 3$ . Then

$$\delta_n(3) = \delta_{n-1}(3) - \frac{1}{n}(\zeta(2) - H_n^{(2)})$$

$$= \ldots = \delta_0(3) - \zeta(2)(1 + \frac{1}{2} + \ldots + \frac{1}{n}) + \frac{1}{n}H_n^{(2)} + \frac{1}{n-1}H_{n-1}^{(2)} + \ldots + H_1^{(2)},$$

i.e.

$$\delta_n(3) = \zeta(3) - \zeta(2)H_n + \sum_{k=1}^{n} \frac{1}{k}H_k^{(2)} . \tag{4.5}$$

## 5. CONVOLUTIONS OF HARMONIC NUMBERS AND STIRLING NUMBERS

In this section we evaluate the series $\tau_n(p)$ defined in (1.5). For every two integers $p \geq 2$ and $k \geq p$ set

$$W_p(k) = \sum_{k_1+\ldots+k_p=k} H_{k_1} H_{k_2} \ldots H_{k_p}, \tag{5.1}$$

where $H_{k_j}$ are the harmonic numbers and $1 \leq k_j \leq k$. The convolutions $W_p(k)$ are the coefficients in the expansion (see (3.3))

$$\left( \frac{-\ln(1-x)}{1-x} \right)^p = \left( \sum_{n=1}^{\infty} H_n x^n \right)^p = \sum_{k=p}^{\infty} W_p(k) x^k . \tag{5.2}$$

Thus, when $p = 2$,

$$W_2(2) = H_1^2 = 1, \; W_2(3) = \sum_{k_1+k_2=3} H_{k_1} H_{k_2} = H_1 H_2 + H_2 H_1 = 2H_2 = 3,$$

$$W_2(k) = \sum_{k_1+k_2=k} H_{k_1} H_{k_2} = \sum_{j=1}^{k-1} H_j H_{k-j} . \tag{5.3}$$

**Proposition 5.1** *For every* $p \geq 2$ *and* $n > p$*,*

$$\tau_n(p) = p!(\zeta(p+1) - H_{n-p}^{(p+1)}). \tag{5.4}$$

*Also,*

$$\tau_p(p) = p!\zeta(p+1). \tag{5.5}$$

*In particular, when* $n = p = 2$ *we have the representation*

$$\zeta(3) \equiv \sum_{k=2}^{\infty} \frac{W_2(k)}{k(k+1)(k+2)} . \tag{5.6}$$

*Proof.* We multiply (3.1) by $W_p(k)$ and sum for $k \geq p$ to obtain (according to (5.2))

$$\tau_n(p) = \int_0^1 \{ \sum_{k=p}^{\infty} W_p(k) x^{k-1} \} (1-x)^n dx = \int_0^1 \left( \frac{-\ln(1-x)}{1-x} \right)^p (1-x)^n \frac{dx}{x}$$

$$= \int_0^1 (-\ln(1-x))^p (1-x)^{n-p} \frac{dx}{x} \ . \tag{5.7}$$

The last integral with the substitution $t = -\ln(1-x)$ turns into

$$\tau_n(p) = \int_0^\infty \frac{t^p e^{-(n-p+1)t}}{1-e^{-t}} dt = p!\,\zeta(p+1, n-p+1), \tag{5.8}$$

where

$$\zeta(s,a) = \sum_{k=0}^\infty \frac{1}{(k+a)^s} = \frac{1}{\Gamma(s)} \int_0^\infty \frac{t^{s-1} e^{-at}}{1-e^{-t}} dt \tag{5.9}$$

(for $Re\, s > 1$, $a > 0$) is the Hurwitz zeta function ([7]). Obviously,

$$\zeta(p+1, n-p+1) = \zeta(p+1) - \sum_{j=1}^{n-p} \frac{1}{j^{p+1}} \tag{5.10}$$

which finishes the proof.

Now we evaluate the last series $\rho_n(m)$ from (1.6).

**Proposition 5.2**. *For every* $n \geq 0$ *and* $m \geq 1$,

$$\rho_n(m) = \zeta(m+1, n+1) = \zeta(m+1) - H_n^{(m+1)}. \tag{5.11}$$

For the proof we use the generating function (see [4, p. 351])

$$[-\ln(1-x)]^m = m! \sum_{k=m}^\infty \begin{bmatrix} k \\ m \end{bmatrix} \frac{x^k}{k!} \ . \tag{5.12}$$

Thus

$$m!\, n! \sum_{k=m}^\infty \frac{1}{k(k+1)\ldots(k+n)} \begin{bmatrix} k \\ m \end{bmatrix} \frac{1}{k!} = \int_0^1 (-\ln(1-x))^m (1-x)^n \frac{dx}{x} \ ,$$

and we continue as in the proof of the previous proposition. The rest is left to the reader.

**Remark 5.3** . The representation

$$\zeta(m+1,n+1) = n!\sum_{k=m}^{\infty}\frac{1}{k(k+1)...(k+n)}\begin{bmatrix}k\\m\end{bmatrix}\frac{1}{k!} \tag{5.13}$$

extend the remarkable formula

$$\zeta(m+1) = \sum_{k=m}^{\infty}\begin{bmatrix}k\\m\end{bmatrix}\frac{1}{k!\,k}, \tag{5.14}$$

as (5.13) turns into (5.14) when $n = 0$. Apparently, (5.14) was first obtained by Jordán Károly (Charles Jordan) [5, p.166, 195]. It has been rediscovered several times after that. Jordán Károly obtained also two similar extensions - see [5, p. 343] and also [7, p. 76] .

Harmonic numbers are connected to Stirling numbers of the first kind as shown, for instance, in [4] and [7, p.57]. This relationship is due to the similarity of their generating functions (3.3), (5.2) and (5.12). One such connection is listed below.

**Corollary 5.4** *From the above two propositions we conclude that for* $n \ge p$ *we have the identity* $\tau_n(p) = p!\,\rho_{n-p}(p)$, *or*

$$\binom{n}{p}\sum_{k=p}^{\infty}\frac{W_p(k)}{k(k+1)...(k+n)} = \sum_{k=p}^{\infty}\frac{1}{k(k+1)...(k+n-p)}\begin{bmatrix}k\\p\end{bmatrix}\frac{1}{k!}. \tag{5.15}$$